\documentclass[preprint,12pt]{elsarticle}

\newtheorem{theorem}{\textbf{Theorem}}

\newtheorem{lemma}{\textbf{Lemma}}

\usepackage{mathrsfs}

\def\a {\alpha}

\def\Z {\mathbb{Z}}
\def\Q {\mathbb{Q}}

\usepackage{amsfonts}
\usepackage{amssymb}

\def\e {\epsilon}
\def\p {\alpha}
\def\d {\delta}

\def\dfrac {\displaystyle\frac}
\def\g {\gamma}

\def\tareesidedbox#1{\setbox0=\hbox{$#1$}\dimen0=\wd0 \advance\dimen0 by3pt\rlap{\hbox{\vrule height8pt width.4pt
 depth2pt \kern-.4pt\vrule height8.4pt width\dimen0 depth-8pt\kern-.4pt \vrule height8pt width.4pt depth2pt}}
 \relax \hbox to\dimen0{\hss$#1$\hss}}

\usepackage{amssymb}
\journal{...}

\begin{document}

\begin{frontmatter}



\title{Powers of two as sum of two generalized Fibonacci numbers}

%\address{Departament of Mathematics, University of British Columbia, Vancouver, Canada}

%\author{Ana Paula Chaves\corref{cor1}}$^{a,}${\ead{apchaves@mat.unb.br}}
%\author{Diego Marques\corref{cor2}$^{a,}$}
%\author{Alain Togb\' e\corref{cor3}$^{b,}$}

\author{Diego Marques\fnref{dmm}}{\ead{diego@mat.unb.br}}
%\author[pt]{Pavel Trojovsk\' y\fnref{ptt}\corref{CA}}{\ead{pavel.trojovsky@uhk.cz}}
\address{Departamento de Matem\' atica, Universidade  de Bras\' ilia, Bras\' ilia, 70910-900, Brazil}
%\address[pt]{Department of Mathematics, University of Hradec Kr\'alov\'e, Faculty of Science, Rokitansk\'eho 62, Hradec Kr\'alov\'e, 500 03, Czech Republic}

%\author{Ana Paula Chaves\corref{cor1}}$^{a,}${\ead{apchaves@mat.unb.br}}
%\author{Diego Marques\corref{cor2}$^{a,}$}
%\author{Alain Togb\' e\corref{cor3}$^{b,}$}

\fntext[dmm]{Supported by FAP-DF and CNPq-Brazil}

\begin{abstract}
For $k\geq 2$, the $k$-generalized Fibonacci sequence $(F_n^{(k)})_{n}$ is defined by the initial values $0,0,\ldots,0,1$ ($k$ terms) and such that each term afterwards is the sum of the $k$ preceding terms. In this paper, we search for powers of two of the form $F_n^{(k)}+F_m^{(k)}$. This work is related to a recent result by Bravo and Luca \cite{Lucanew} concerning the case $k=2$.
\end{abstract}


\end{frontmatter}

% \linenumbers


%&&&&&&&&&&&&&&&&&&&&&&&&&&&&&&&&&&&&&&&&&&&&&&&&&&&&&&&&&&&&&%
\section{Introduction}
%&&&&&&&&&&&&&&&&&&&&&&&&&&&&&&&&&&&&&&&&&&&&&&&&&&&&&&&&&&&&&%

Let $(F_n)_{n\geq 0}$ be the Fibonacci sequence given by $F_{n+2}=F_{n+1}+F_n$, for $n\geq 0$, where $F_0=0$ and $F_1=1$. These numbers are well-known for possessing amazing properties (consult \cite{FIB3} together with its very extensive annotated bibliography for additional references and history). 

The problem of searching for Fibonacci numbers of a particular form has a very rich history, see for example \cite{bugeaud} and references therein. In this same paper, Bugeaud et al \cite[Theorem 1]{bugeaud} showed that $0,1,8,144$ are the only perfect powers in the Fibonacci sequence. Other related papers searched for Fibonacci numbers of the forms $px^2+1,px^3+1$ \cite{px}, $k^2+k+2$ \cite{k}, $p^a\pm p^b+1$ \cite{p+1}, $p^a\pm p^b$ \cite{p}, $y^t\pm 1$ \cite{yt}, $q^ky^t$ \cite{qy} and $2^a+3^b+5^c$ \cite{MT}. 

Let $k\geq 2$ and denote $F^{(k)}:=(F_n^{(k)})_{n\geq -(k-2)}$, the \textit{$k$-generalized Fibonacci sequence} whose terms satisfy the recurrence relation
\begin{equation}\label{rec}
F_{n+k}^{(k)}=F_{n+k-1}^{(k)}+F_{n+k-2}^{(k)}+\cdots + F_{n}^{(k)},
\end{equation}
with initial conditions $0,0,\ldots,0,1$ ($k$ terms) and such that the first nonzero term is $F_1^{(k)}=1$.

The above sequence is one among the several generalizations of Fibonacci numbers. Such a sequence is also called $k$-\textit{step Fibonacci sequence}, \textit{Fibonacci $k$-sequence}, or $k$-\textit{bonacci sequence}. Clearly for $k=2$, we obtain the classical Fibonacci numbers $(F_n)_n$, for $k=3$, the \textit{Tribonacci} numbers $(T_n)_n$, for $k=4$, the \textit{Tetranacci} numbers $(Q_n)_n$, etc.

Recently, these sequences have been the main subject of many works. We refer to \cite{BLP} for results on the largest prime factor of $F_n^{(k)}$ and we refer to \cite{BL2} for the solution of the problem of finding powers of two belonging to these sequences. In 2013, two conjectures concerning these numbers were proved. The first one, proved by Bravo and Luca \cite{BL} is related to \textit{repdigits} (i.e., numbers with only one distinct digit in its decimal expansion) among $k$-bonacci numbers (proposed by Marques \cite{util}) and the second one, a conjecture (proposed by Noe and Post \cite{noe}) about coincidences between terms of these sequences, proved independently by Bravo-Luca \cite{BL0} and Marques \cite{Marques0} (see \cite{spa} for results on the spacing between terms of these sequences). Also, Bravo and Luca \cite{BLnew} found all repdigits which are sum of two $k$-bonacci numbers. Recently, all generalized Fibonacci numbers of the form $2^a+3^b+5^c$ (with $\max\{a,b\}\leq c$) and $\sum_{j=1}^k(2^j)^{n_j}$ (with $\max_{1\leq i\leq k-1}\{n_i\}\leq n_k$) were found, see \cite{Mbr} and \cite{ruiz}, respectively.

Very recently, Bravo and Luca \cite{Lucanew} found all solutions of $F_n+F_m=2^t$ (see also \cite{BLJIS}).The aim of this paper is to study their equation in the $k$-bonacci context. More precisely, our main results are the following
\begin{theorem}\label{main}
The Diophantine equation
\begin{equation}\label{F}
F_n^{(k)}+F_m^{(k)}=2^t
\end{equation}
has no solution $(n,m,t,k)$, with $2\leq m<n\neq t+2$. 

\end{theorem}

We remark that the condition $n> m$ is to avoid (by symmetry) the case $n=m$ which leads to the equation $F_n^{(k)}=2^{t-1}$ already solved in \cite{BL2}. 

Note that if (\ref{F}) has solution with $2\leq m<n$, then $n=t+2$ (in fact, there are infinitely many families of solutions). Our next result treats this case (as usual, $[a,b]$ denotes $\{a,a+1,\ldots,b\}$, for integers $a<b$).

\begin{theorem}\label{main2}
There is no integer solution $(n,m,t,k)$ of Eq. (\ref{F}) with $2\leq m<n=t+2$, when $(n,m)\in [2,k+1]^2$ or $[k+2,2k+2]^2$. In the case of $(n,m)\in [k+2,2k+2]\times [2,k+1]$ we have that for any given $s\geq 1$, then
\[
(n,m,t,k)=(2^s+k,2^s+s-1,2^s+k-2,k)
\]
is solution, for all $k\geq 2^s+s-2$.
\end{theorem}

In fact, when $(n,m)\in [k+2,2k+2]\times [2,k+1]$,  we have that $(n,m,t,k)=(t+2,m,t,k)$ is a solution only if it satisfies
\[
t-k+2=2^{m+k-t-1}.
\]

Finding all the solutions of the above equation (and so all solutions for (\ref{F})) seems to be an untreatable problem. In fact, even in the particular case $t=2p-1$, $m=p+1$, $p< k$ with $2^p-1$ prime, we must find all solutions of the Diophantine equation $2p-k+1=2^{k-p+1}$ (which is related to even perfect numbers among generalized Fibonacci numbers, see \cite{Dnew}). Clearly, there is no solution when $k\not\equiv 3\pmod 4$, but when $k\equiv 3\pmod 4$, the possible solutions of the equation are related to Mersenne primes of the form $2^{2^{\ell}+\ell-1}-1$ and it seems out of reach to find all these primes, if any (there is no  such primes when $\ell<5\cdot 10^5$).

Our proof of Theorem \ref{main} combines lower bounds for linear forms in three logarithms, a variation of a Dujella and Peth\H  o reduction lemma and a fruitful method (doubly used) developed by Bravo and Luca concerning approximation of some convenient number (related to the dominant root of the characteristic polynomial of $F_n^{(k)}$) by a power of $2$. The proof of Theorem \ref{main2} uses closed formulas for $F_j^{(k)}$. The calculations in this paper took roughly $15$ days on $2.5$ GHz Intel Core i5 4GB Mac OSX.

%&&&&&&&&&&&&&&&&&&&&&&&&&&&&&&&&&&&&&&&&&&&&&&&&&&&&&&&&&&&&&%
\section{Auxiliary results}\label{sec2}
%&&&&&&&&&&&&&&&&&&&&&&&&&&&&&&&&&&&&&&&&&&&&&&&&&&&&&&&&&&&&&%

Before proceeding further, we shall recall some facts and tools which will be used
after.

We know that the characteristic polynomial of $(F_n^{(k)})_n$ is 
\[
\psi_k(x):=x^k-x^{k-1}-\cdots -x-1
\]
and it is irreducible over $\Q[x]$ with just one zero outside the unit circle. That single zero is located between $2(1-2^{-k})$ and $2$ (as can be seen in \cite{wolf}). Also, in a recent paper, Dresden and Du \cite[Theorem 1]{dres} gave a simplified ``Binet-like"\ formula for $F_n^{(k)}$:
\begin{equation}\label{binet}
F_n^{(k)}=\displaystyle\sum_{i=1}^k\dfrac{\alpha_i-1}{2+(k+1)(\a_i-2)}\a_i^{n-1},
\end{equation}
for $\alpha=\alpha_1,\ldots,\a_k$ being the roots of $\psi_k(x)$. Also, it was proved in \cite[Lemma 1]{BL} that
\begin{equation}\label{Lu}
\alpha^{n-2}\leq F_n^{(k)}\leq \alpha^{n-1},\ \mbox{for\ all}\ n\geq 1,
\end{equation}
where $\a$ is the dominant root of $\psi_k(x)$. Moreover, Bravo and Luca \cite[Lemma 1]{BL2} proved that $F_n^{(k)}\leq 2^{n-2}$, for all $n\geq 2$. Also, the contribution of the roots inside the unit circle in formula (\ref{binet}) is almost trivial. More precisely, it was proved in \cite{dres} that
\begin{equation}\label{small}
|F_n^{(k)}-g(\alpha,k)\alpha^{n-1}|<\dfrac{1}{2},
\end{equation}
where we adopt throughout the notation $g(x,y):=(x-1)/(2+(y+1)(x-2))$.

Another tool to prove our theorem is a lower bound for a linear form logarithms {\it \`a la Baker} and such a bound was given by the following result of Matveev (see \cite{matveev} or Theorem 9.4 in \cite{bugeaud}).
%\begin{lemma}\label{lemma0}(Baker-W$\ddot{\mbox{u}}$stholz)
%Let
%$$\Lambda=b_1\log \alpha_1+\cdots + b_s\log \alpha_s,$$
%where $b_i\in \mathbb{Z}$ and $H(\alpha)$ (so-called the \emph{height} of $\alpha$) is the maximum of absolute values of coefficients of the minimal polynomial of $\alpha$ (over $\mathbb{Z}$). Suppose that $d=[\mathbb{Q}(\alpha_1,\ldots,\alpha_s):\mathbb{Q}]$, $A_i\geq \max\{H(\alpha_i),e\}$ for $i=1,\ldots,s$ and $B\geq \max\{|b_1|,\ldots,|b_s|,e\}$. If $\Lambda\neq 0$, then
%$$|\Lambda|>\exp (-(16ds)^{2s+4}\cdot \log A_1\cdots \log A_s\cdot \log B)$$
%\end{lemma}

%To prove Theorem \ref{thm:main2}, Binet's formulas will allow us to obtain linear forms in three logarithms and then we will determine lower bounds for these linear forms. The bounds available for linear forms in three logarithms are substantially better than those available for general linear forms in logarithms. From the main result of Matveev \cite{matveev}, we deduce the following lemma.
\begin{lemma}\label{lemma1}
Let $\gamma_1,\ldots,\gamma_t$ be real algebraic numbers and let $b_1,\ldots,b_t$ be nonzero rational integer numbers. Let $D$ be the degree of the number field $\Q(\gamma_1,\ldots,\gamma_t)$ over $\Q$ and let $A_j$ be a real number satisfying
\begin{center}
$A_j\geq \max\{Dh(\g_j),|\log \g_j|,0.16\}$, for $j=1,\ldots,t$.
\end{center}
Assume that
$$
B\geq \max\{|b_1|,\ldots,|b_t|\}.
$$
If $\g_1^{b_1}\cdots \g_t^{b_t}\neq 1$, then
$$
|\g_1^{b_1}\cdots \g_t^{b_t}-1|\geq \exp(-1.4\cdot 30^{t+3}\cdot t^{4.5}\cdot D^2(1+\log D)(1+\log B)A_1\cdots A_t).
$$
\end{lemma}
As usual, in the previous statement, the \textit{logarithmic height} of an $s$-degree algebraic number $\g$ is defined as
$$
h(\g)=\dfrac{1}{s}(\log |a|+\displaystyle\sum_{j=1}^s\log \max\{1,|\g^{(j)}|\}),
$$
where $a$ is the leading coefficient of the minimal polynomial of $\gamma$ (over $\mathbb{Z}$) and $(\g^{(j)})_{1\leq j\leq s}$ are the conjugates of $\g$ (over $\Q$).
% It can be deduced from Lemma 1 in \cite{MT2}, but we shall prove it for the convenience of the reader.

After finding an upper bound on $n$ which is too large for practical purposes, the next step is to reduce it. For that, our last ingredient can be found in \cite[Lemma 4]{BL} and it is a variant of the famous Dujella and Peth\H{o} \cite[Lemma 5 (a)]{dujella} reduction lemma. For a real number $x$, we use  $\parallel x \parallel= \min\{|x-n|:n\in \Z\}$ for the distance from $x$ to the nearest integer.
\begin{lemma}\label{lemma2}
Suppose that $M$ is a positive integer. Let $p/q$ be a convergent of the continued fraction expansion of the irrational number $\gamma$ such that $q > 6M$ and let $A, B$ be some real numbers with $A>0$ and $B>1$. Let $\epsilon=\parallel \mu q \parallel-M\parallel \gamma q \parallel$, where $\mu$ is a real number. If $\epsilon>0$, then there is no solution to the inequality
$$
0<m\gamma -n+\mu < A\cdot B^{-k}
$$
in positive integers $m,n$ and $k$ with
\begin{center}
$m\leq M$ and $k\geq \dfrac{\log(Aq/\epsilon)}{\log B}$.
\end{center}
\end{lemma}

Now, we are ready to deal with the proof of the theorems.

%&&&&&&&&&&&&&&&&&&&&&&&&&&&&&&&&&&&&&&&&&&&&&&&&&&&&&&&&&&&&&%
\section{The proof of the theorems} \label{sec3}
%&&&&&&&&&&&&&&&&&&&&&&&&&&&&&&&&&&&&&&&&&&&&&&&&&&&&&&&&&&&&&%

\subsection{The proof of Theorem \ref{main}}

In light of \cite{Lucanew}, throughout this paper we shall assume $k\geq 3$.

%---------------------------------------------------------------------%
\subsubsection{Upper bounds for $n, m$ and $t$ in terms of $k$}\label{sec31}
%---------------------------------------------------------------------% 

In this section, we shall prove the following result
\begin{lemma}\label{t1}
If $(n,m,t,k)$ is an integer solution of Diophantine equation (\ref{F}), with $n>m$, then
\begin{equation}\label{c<k}
\max\{m,t\}<n<1.2\cdot 10^{27}k^7\log^5 k.
\end{equation}
\end{lemma}

\noindent
{\bf Proof.}  We use Eq. (\ref{F}) together with (\ref{binet}) to obtain
\begin{equation}\label{s}
g(\a,k)\a^{n-1}-2^t=-E_k(n)-F_m^{(k)}<0,
\end{equation}
where $E_k(n):=\sum_{i=2}^kg(\a_i,k)\a_i^{n-1}$. Thus, by (\ref{Lu}) and the fact that $|E_k(n)|<1/2$, we have
\[
|g(\a,k)\a^{n-1}-2^t|<\frac{1}{2}+\a^{m-1}<2\a^{m-1}
\]
and then
\begin{equation}\label{Est1}
\left|\dfrac{2^{t}}{g(\a,k)\a^{n-1}}-1\right|<\dfrac{4}{\a^{n-m}},
\end{equation}
where we used that $g(\a,k)>1/\a>1/2$ (see \cite[Section 2]{CH}).

In order to use Lemma \ref{lemma1}, we take $t:=3$,
\[
\g_1:=2,\ \g_2:=\a,\ \g_3:=g(\a,k)
\]
and 
\[
b_1:=t,\ b_2:= -(n-1),\ b_3:= -1.
\]
For this choice, we have $D=[\Q(\a):\Q]= k$. In \cite[p. 73]{BL2}, an estimate for $h(g(\a,k))$ was given. More precisely, it was proved that
\[
h(\g_3)=h(g(\a,k)) < \log (4k+4).
\]
Note that $h(\g_1)=\log 2$ and $h(\g_2)< 0.7/k$. Thus, we can take $A_1:= k\log 2, A_2:= 0.7$ and $A_3:=k\log (4k+4)$.  

Note that $\max\{|b_1|,|b_2|,|b_3|\}=\max\{t,n-1\}$. By using the inequalities in (\ref{Lu}), we get $2^t=F_n^{(k)}+F_m^{(k)}< 2F_n^{(k)}<2^n$. Thus, we can choose $B:=n-1$. Since $2^{t}g(\a,k)\a^{n-1}-1>0$ (by (\ref{s})), we are in position to apply Lemma \ref{lemma1}. This lemma together with a straightforward calculation gives
\begin{equation}\label{Est2}
\left|\dfrac{2^{t}}{g(\a,k)\a^{n-1}}-1\right| >  \exp(-7.3\cdot 10^{11}k^{4}\log n\log^2 k),
\end{equation}
where we used that $1+\log k<2\log k$, for $k\geq 3$, $1+\log (n-1)<2\log n$, for $c\geq 2$, and $\log (4k+4)<2.6\log k$, for $k\geq 3$.

By combining (\ref{Est1}) and (\ref{Est2}), we obtain
\begin{equation}\label{n-m}
n-m<1.5\cdot 10^{12}k^{4}\log^2 k\log n.
\end{equation}

Now, we shall see Eq. (\ref{F}) by another viewpoint
\[
g(\a,k)\a^{n-1}+g\a^{m-1}-2^t=-E_k(n)-E_k(m).
\]
After some straightforward manipulation, we arrive at
\begin{equation}\label{t1}
\left|1-\dfrac{2^t\a^{-(n-1)}}{g(\a,k)(1+\a^{m-n})}\right|<\dfrac{2}{(1.3)^n},
\end{equation}
where we used that $(1.3)^{n}<\a^{n-1}$, for $k\geq 3$ and $n>1$.

In order to use Lemma \ref{lemma1}, we take $t:=3$,
\[
\g_1:=2,\ \g_2:=\a,\ \g_3:=g(\a,k)(1+\a^{m-n})
\]
and 
\[
b_1:=t,\ b_2:= -(n-1),\ b_3:= -1.
\]
For this choice, we have $D=[\Q(\a):\Q]= k$. We use that $h(xy)\leq h(x)+h(y), h(x+y)\leq h(x)+h(y)+\log 2$ and $h(x^s)=|s|h(x)$ (for $s\in \Z$) to obtain
\[
h(\g_3)<\log (8k+8)+(n-m)\frac{\log 2}{k}.
\]
Thus, we can take $A_1:= k\log 2, A_2:= 0.7$ and $A_3:=k\log (8k+8)+(n-m)\log 2$.  Now, we need to prove that the left-hand side of (\ref{t1}) is not zero. Suppose the contrary, then we have that $2^t=g(\a,k)(\a^{n-1}+\a^{m-1})$. By conjugating the previous relation in the splitting field of $\psi_k(x)$, we obtain $2^{t}=\a_i^{m-1}g(\a_i,k)(\a_i^{n-1}+\a_i^{m-1})$, for $i=1,\ldots,k$. However, when $i>1$, $|\a_i|<1$ and $|g(\a_i,k)|\leq 1$. But this leads to the following absurdity
\[
2^{t}=|\a_i^{m-1}g(\a_i,k)(\a_i^{n-m}+1)|\leq |\a_i^{m-1}||g(\a_i,k)|(|\a_i|^{n-1}+|\a_i|^{m-1})|\leq 2,
\]
since $t>1$. Then we are in position to apply Lemma \ref{lemma1}. Therefore
\[
\left|1-\dfrac{2^t\a^{-(n-1)}}{g(\a,k)(1+\a^{m-n})}\right|>\exp(-2.8\cdot 10^{11}k^3\log k\log n (k\log (8k+8)+(n-m)\log 2)).
\]
Since $\log (8k+8)\leq 3.4\log k$, we have $k\log (8k+8)+(n-m)\log 2\leq 3.4(n-m+k)\log k$. On the other hand, by (\ref{n-m}), $n-m+k<1.6\cdot 10^{12}k^{4}\log^2 k\log n$. Therefore
\[
\left|1-\dfrac{2^t\a^{-(n-1)}}{g(\a,k)(1+\a^{m-n})}\right|>\exp(-4.5\cdot 10^{23}k^7\log^3 k\log^2 n).
\]
By combining the previous inequality with (\ref{t1}), we get
\[
\dfrac{n}{\log^2 n}<1.8\cdot 10^{23}k^7\log^3 k.
\]
Since the function $x\mapsto x/\log^2 x$ is increasing for $x>e$, then it is a simple matter to prove that
\begin{equation}\label{key}
\dfrac{x}{\log^2 x}<A\ \ \mbox{implies\ that}\ \ x<2A\log^2 A\  (\mbox{for}\ A\geq 2252750).
\end{equation}
In fact, suppose the contrary, i.e., $x\geq 2A\log^2 A$. Then
\[
\dfrac{x}{\log^2 x}\geq \dfrac{2A\log^2 A}{\log^2 (2A\log^2 A)}>A,
\]
which contradicts our inequality. Here we used that $\log^2 (2A\log^2 A)<2\log^2 A$, for $A\geq 2252750$.

Thus, using (\ref{key}) for $x:=n$ and $A:=1.8\cdot 10^{23}k^7\log^3 k$, we have that
\[
m<2(1.8\cdot 10^{23}k^7\log^3 k)\log^2 (1.8\cdot 10^{23}k^7\log^3 k).
\]
A straightforward calculation gives
\[
n<1.2\cdot 10^{27}k^7\log^5 k,
\]
where we used that $\log(1.8\cdot 10^{23}k^7\log^3 k)<56\log k$.
\qed

\subsubsection{The small cases: $3\leq k\leq 321$}\label{small2}

In this section, we shall prove the following result
\begin{lemma}\label{t2}
If $(n,m,t,k)$ is an integer solution of Diophantine equation (\ref{F}), with $3\leq k\leq 321$ and $2\leq m<n\neq t+2$.
\end{lemma}
\noindent
{\bf Proof.} By using (\ref{s}) and (\ref{Est1}), we have that
\[
0<t\log 2-(n-1)\log \a+\log (1/g(\a,k))<4\a^{-(n-m)}.
\]
Dividing by $\log \a$, we obtain
\begin{equation}\label{DP}
0<t\g_k-(n-1)+\mu_k<7.2\cdot \a^{-(n-m)},
\end{equation}
where $\g_k=\log 2/\log \a^{(k)}$ and $\mu_k=\log (1/g(\a^{(k)},k))/\log \a^{(k)}$. Here, we added the superscript to $\a$ for emphasizing its dependence on $k$.

We claim that $\g_k$ is irrational, for any integer $k\geq 2$. In fact, if $\g_k=p/q$, for some positive integers $p$ and $q$, we have that $2^q=(\a^{(k)})^p$ and we can conjugate this relation by some automorphism of the Galois group of the splitting field of $\psi_k(x)$ over $\Q$ to get $2^q=|(\a_i^{(k)})^p|<1$, for $i>1$, which is an absurdity, since $q\geq 1$. Let $q_{\ell,k}$ be the denominator of the $\ell$-th convergent of the continued fraction of $\g_k$. Taking $M_k:=1.2\cdot 10^{27}k^7\log^5 k$, we use \textit{Mathematica} \cite{math} to get
\[
\displaystyle\min_{3\leq k\leq 321}q_{120,k}>7.7\cdot 10^{51}>6M_{321}.
\]
Also
\[
\displaystyle\max_{3\leq k\leq 321}q_{120,k}<4.2\cdot 10^{162}.
\]

Define $\e_k:=\parallel \mu_kq_{120,k}\parallel-M_k\parallel \g_kq_{120,k} \parallel$, for $3\leq k\leq 321$, we get (again using Mathematica)
\[
\displaystyle\min_{3\leq k\leq 321}\e_{k}>2.3\cdot 10^{-37}.
\]

Note that the conditions to apply Lemma \ref{lemma2} are fulfilled for $A=7.2$ and $B=\a$, and hence there is no solution to inequality (\ref{DP}) (and then no solution to the Diophantine equation (\ref{F})) for $t$ and $n-m$ satisfying
\begin{center}
$t<M_k$ and $n-m \geq \dfrac{\log (Aq_{120,k}/\e_k)}{\log B}$.
\end{center}
Since $t<M_k$ (Lemma \ref{t1}), then
\[
n-m < \dfrac{\log (Aq_{120,k}/\e_k)}{\log B}\leq \dfrac{\log (7.2\cdot 4.2\cdot 10^{162}/2.3\cdot 10^{-37})}{\log \a}\leq 843.978\ldots.
\]

Therefore $n-m\leq 843.$ Now, set $\Gamma=t\log 2-(n-1)\log \a+\log \phi(k,n-m)$, where $\phi(k,t)=1/(g(\a,k)(1+\a^{-t}))$. We know that $\Gamma\neq 0$ (by a previous argument), so let us suppose that $\Gamma>0$ (the case $\Gamma<0$ can be handled in much the same way). Thus $\Gamma<e^{\Gamma}-1<2/(1.3)^n$ (where we used (\ref{t1})). Therefore
\begin{equation}\label{DP2}
0<t\g_k-(n-1)\log \a+\mu^*_{k,n-m}<2\cdot (1.3)^{-n},
\end{equation}
where $\mu^*_{k,n-m}=\log \phi(k,n-m)/\log \a^{(k)}$.

Define $\e^*_{k,\ell}:=\parallel \mu^*_{k,\ell}q_{120,k}\parallel-M_k\parallel \g_kq_{120,k} \parallel$, for $3\leq k\leq 321$ and $1\leq \ell\leq 843$, we get (again using Mathematica)
\[
\displaystyle\min_{k\in [3,321], \ell\in [1,843]}\e^*_{k,\ell}>5.6\cdot 10^{-91}.
\]

We apply again Lemma \ref{lemma2} for $A=2$ and $B=1.3$, and hence there is no solution to inequality (\ref{DP2}) (and then no solution to Eq. (\ref{F})) for $t$ and $n$ satisfying
\begin{center}
$t<M_k$ and $n\geq \dfrac{\log (Aq_{120,k}/\e^*_{k,\ell})}{\log B}$.
\end{center}
Since $t<M_k$ (Lemma \ref{t1}), then
\[
n < \dfrac{\log (Aq_{120,k}/\e^*_{k,\ell})}{\log B}\leq \dfrac{\log (2\cdot 4.2\cdot 10^{162}/5.6\cdot 10^{-91})}{\log (1.3)}=2265.83\ldots.
\]

Therefore, we have $\max\{m,t\}<n\leq 2265$, $n\neq t+2$ and $3\leq k\leq 321$. Now, we use a Mathematica routine which does not return any solution for (\ref{F}) with the previous conditions. \qed

\subsubsection{The proof}

In this section, we shall complete the proof of Theorem \ref{main} by proving the following result
\begin{lemma}\label{t3}
There is no integer solution $(n,m,t,k)$ of Eq. (\ref{F}), with $k\geq 322$ and $2\leq m<n\neq t+2$. 
\end{lemma}
\noindent
{\bf Proof.} By using that $F_n^{(k)}\leq 2^{n-2}$, we get $2^t=F_n^{(k)}+F_m^{(k)}<2^{n-1}$ and so $n>t+1$. Hence, we have $n>t+2$ (since $n\neq t+2$).

Since $k\geq 322$, we have
\begin{equation}\label{exact}
n<1.2\cdot 10^{27}k^7\log^5 k<2^{k/2}.
\end{equation}

Now, we use doubly a key argument due to Bravo and Luca \cite[p. 77-78]{BL2}. Let $\d_j=\a^{j-1}-2^{j-1}$, for $j\in \{n,m\}$, and $\eta=g(\a,k)-1/2$. Then, Bravo and Luca proved that $|\d_j|<2^j/2^{k/2}$ and $|\eta|<2k/2^k$. Thus, we deduce that
\[
2^{n-2}=g(\a,k)\a^{n-1}-2^{n-1}\eta - \dfrac{\d_n}{2} - \d_n \eta
\]
and
\[
2^{m-2}=g(\a,k)\p^{m-1}-2^{m-1}\eta - \dfrac{\d_m}{2} - \d_m \eta.
\]
If $n\leq k+1$, then $F_n^{(k)}+F_m^{(k)}=2^{n-2}+2^{m-2}$ cannot be a power of two (since $m<n$). So, we may suppose $n>k+1$. After some manipulations, we have
\[
|2^{n-2}+2^{m-1}-2^t|<\dfrac{5\cdot 2^{n-2}}{2^{k/2}}
\]
and so
\[
|1-\dfrac{1}{2^{n-t-2}(1+2^{m-n})}|<\dfrac{5}{2^{k/2}}.
\]
However, $2^{n-t-2}(1+2^{m-n})>2^{n-t-2}\geq 2$, since $n>t+2$. Thus
\[
\dfrac{5}{2^{k/2}}>|1-\dfrac{1}{2^{n-t-2}(1+2^{m-n})}|=1-\dfrac{1}{2^{n-t-2}(1+2^{m-n})}>1/2
\]
yielding $2^{k/2}<10$ which is false for $k\geq 322$.
\qed

\subsection{The proof of Theorem \ref{main2}}

To deal with those cases, we shall use the facts that $F_i^{(k)}=2^{i-2}$, if $i\in [2,k+1]$ and $F_{j}^{(k)}=2^{j-2}-(j-k)2^{j-k-3}$, for $j\in [k+2,2k+2]$ (see \cite[Theorem 2.2]{coo}).

First, note that there is no solution when $2\leq m<n\leq k+1$ (so $F_m^{(k)}=2^{m-2}$ and $F_n^{(k)}=2^{n-2}$), because the sum of two distinct powers of two cannot be a power of two.

Also, there is no solution when $k+2\leq m,n\leq 2k+2$ and $n=t+2$. In fact, it is easy to see this when $k\in [3,6]$. So, we may assume that $k\geq 7$. If $(n,m,t,k)$ is a solution of (\ref{F}) (with $n=t+2$), then we use that $F_{j}^{(k)}=2^{j-2}-(j-k)2^{j-k-3}$, for $j\in \{m,n\}$, to arrive at 
 \[
 2^{k+1}=(n-k)2^{n-m}-(m-k).
 \]
Since $n-m\leq 2k+2-(k+2)=k$, then $2^{n-m}\mid m-k$ and in particular, $2^{n-m}\leq m-k\leq k+2$. On the other hand, $(n-k)2^{n-m}=2^{k+1}+m-k>2^{k+1}$ and also $n-k\leq k+2<2^{k/2}$ (since $k\geq 7$). Summarizing, we have 
\[
2^{n-m+k/2}>(n-k)2^{n-m}>2^{k+1}
\]
yielding $n-m>k/2+1$. We then get to the absurdity that $2^{k/2+1}<2^{n-m}\leq k+2$. 

Now, we shall show that $(n,m,t,k)=(2^s+k,2^s+s-1,2^s+k-2,k)$ is solution, for all $k\geq 2^s+s-2$. In fact, we have that $2\leq m=2^s+s-1\leq k+1$ and $k+2\leq n=2^s+k\leq 2k+2-s<2k+2$ (here we used that $2^s\leq k-s+2$). Thus, we can use the closed formulas to get
\begin{eqnarray*}
F_n^{(k)}+F_m^{(m)} & = & 2^{n-2}-(n-k)2^{n-k-3}+2^{m-2}\\
 & = & 2^{2^s+k-2}-(2^s+k-k)2^{2^s+k-k-3}+2^{2^s+s-3}\\
  & = & 2^{2^s+k-2}=2^t.
\end{eqnarray*}
\qed

%\section*{Acknowledgement}


\begin{thebibliography}{9999}



\bibitem{BL2} J. J. Bravo, F. Luca, Powers of two in generalized Fibonacci sequences, \textit{Rev. Colombiana Mat.} \textbf{46} (2012), 67--79.

\bibitem{BL0} J. J. Bravo, F. Luca, Coincidences in generalized Fibonacci sequences, \textit{J. Number Theory.} \textbf{133} (2013), 2121--2137.

\bibitem{BLP} J. J. Bravo, F. Luca, On the largest prime factor of the $k$-Fibonacci numbers. \textit{Int. J. Number Theory} \textbf{9} (2013), 1351-1366.

\bibitem{BL} J. J. Bravo, F. Luca, On a conjecture about repdigits in $k$-generalized Fibonacci sequences, \textit{Publ. Math. Debrecen} {\bf 82} Fasc. 3-4 (2013).

\bibitem{BLnew} J. J. Bravo, F. Luca, Repdigits as sums of two $k$-Fibonacci numbers. 
\textit{Monatsh. Math.} {\bf 175} (2014), 1-21.

\bibitem{BLJIS} J. J. Bravo, F. Luca, Powers of two as sums of two Lucas numbers. \textit{ J. Integer Seq.} {\bf 17} (2014), Article 14.8.3


\bibitem{Lucanew} J. J. Bravo, F. Luca, On the Diophantine equation $F_n + F_m = 2^a$, preprint.



\bibitem{bugeaud} Y. Bugeaud, M. Mignotte, S. Siksek, Classical and modular approaches to exponential Diophantine equations I. Fibonacci and Lucas powers, \textit{Annals of Math.} \textbf{163} (2006), 969-1018.

\bibitem{yt} Y. Bugeaud, M. Mignotte, F. Luca, S. Siksek, Fibonacci numbers at most one away from a perfect power, \textit{Elem. Math.} \textbf{63} (2008), 65--75.

\bibitem{qy} Y. Bugeaud, M. Mignotte, S. Siksek, Sur les nombres de Fibonacci de la forme $q^ky^t$. \textit{C. R. Math. Acad. Sci. Paris} \textbf{339} (2004), no. 5, 327--330. 



\bibitem{CH} A. P. Chaves, D. Marques, A Diophantine equation related to the sum of squares of consecutive $k$-generalized Fibonacci numbers. \textit{The Fibonacci Quart.} \textbf{52} (2014), no. 1, 70-74.

%\bibitem{dub} H. Dubner. Generalized Cullen numbers. \textit{J. Recreat. Math.}, \textbf{21} (1989), 190-194.

\bibitem{coo} C. Cooper, F. T. Howard. Some identities for $r$-Fibonacci numbers, \textit{Fibonacci Quart.} \textbf{49} (2011), no. 3, 231-243.

\bibitem{dres} G. P. Dresden, Z. Du, A simplified Binet formula for $k$-generalized Fibonacci numbers, \textit{J. Integer Seq. }{\bf 17} (2014) No. 4, Article 14.4.7.



\bibitem{dujella} A. Dujella, A. Peth\H{o}, A generalization of a theorem of Baker and Davenport, \textit{Quart. J. Math. Oxford Ser.} (2) \textbf{49} (1998), 291--306.


\bibitem{Dnew} V. Fac\' o, D. Marques, Even perfect numbers among generalized Fibonacci sequences. To appear in \textit{Rend. Circ. Mat. Palermo}.




%\bibitem{uber} F. Heppner, $\ddot{\mbox{U}}$ber Primzahlen der Form $n2^n+ 1$ bzw. $p2^p+ 1$, \textit{Monatsh. Math.} \textbf{85} (1978), 99-103.



%\bibitem{fib2} D. Kalman, R. Mena, The Fibonacci Numbers--Exposed, \textit{Math. Mag.} \textbf{76}: 3 (2003), 167-181.



\bibitem{FIB3} T. Koshy, \textit{Fibonacci and Lucas Numbers with Applications}, Wiley, New York, 2001.

%\bibitem{Laurent} M.~Laurent, Linear forms in two logarithms and interpolation determinants II, \textit{Acta Arith.} \textbf{133.4} (2008), 325--348.

\bibitem{k} F. Luca, Fibonacci numbers of the form $k^2+k+2$. \textit{Applications of Fibonacci numbers,} Vol. 8 (Rochester, NY, 1998), 241--249, Kluwer Acad. Publ., Dordrecht, 1999.

\bibitem{p+1} F. Luca, L. Szalay, Fibonacci numbers of the form $p^a\pm p^b+1$. \textit{Fibonacci Quart.} \textbf{45} (2007), no. 2, 98-103.

\bibitem{p} F. Luca, P. St\u anic\u a, Fibonacci numbers of the form $p^a\pm p^b$. Proceedings of the Eleventh International Conference on Fibonacci Numbers and their Applications. \textit{Congr. Numer.} \textbf{194} (2009), 177-183.

\bibitem{ruiz} F. Luca, C. A. G. Ru\' iz, $k$-generalized Fibonacci numbers of the form $1+2^{n_1}+4^{n^2}+\cdots + (2^k)^{n_k}$, \textit{Math. Commun.} \textbf{19} (2014), 321-332.


\bibitem{Marques0} D. Marques, The proof of a conjecture concerning the intersection of $k$-generalized Fibonacci sequences, \textit{Bull. Brazilian Math. Soc.} {\bf 44} (3) (2013), 455-468.


\bibitem{spa} D. Marques, On the spacing between terms of generalized Fibonacci sequences, \textit{Coll. Math.} \textbf{134} (2014), 267-280.


\bibitem{util} D. Marques, On $k$-generalized Fibonacci numbers with only one distinct digit, To appear in \textit{Util. Math.}



\bibitem{MT} D. Marques, A. Togb\' e, Fibonacci numbers of the form $2^a+3^b+5^c$,  \textit{Proc. Japan Acad. Ser. A Math. Sci}, \textbf{89} (2013), 47--50.

\bibitem{Mbr} D. Marques, On generalized Fibonacci numbers of the form $2^a+3^b+5^c$. To appear in \textit{Bull. Brazilian Math. Soc}.




\bibitem{matveev} E. M. Matveev, An explicit lower bound for a homogeneous rational linear form in logarithms of algebraic numbers, II, \textit{Izv. Ross. Akad. Nauk Ser. Mat.} {\bf 64} (2000), 125--180. English translation in \textit{Izv. Math.} {\bf 64} (2000), 1217--1269.



\bibitem{noe} T. D. Noe and J. V. Post, Primes in Fibonacci $n$-step and Lucas $n$-step sequences, \textit{J. Integer Seq.}, \textbf{8} (2005), Article 05.4.4.

\bibitem{px} N. Robbins, Fibonacci numbers of the forms $pX^2\pm 1,\;pX^3\pm 1,$ where $p$ is prime. \textit{Applications of Fibonacci numbers} (San Jose, CA, 1986), 77--88, Kluwer Acad. Publ., Dordrecht, 1988.


\bibitem{math} Wolfram Research, Inc., Mathematica, Version {\tt 7.0}, Champaign, IL (2008).

\bibitem{wolf} A. Wolfram, Solving generalized Fibonacci recurrences, \textit{Fibonacci Quart.} {\bf 36} (1998),
129--145. 














\end{thebibliography}
\end{document}